\documentclass[11pt, a4paper, twoside]{amsart}
\usepackage[centering, totalwidth = 370pt, totalheight = 625pt]{geometry}
\usepackage{amssymb, amsmath, amsthm, enumerate, microtype, stmaryrd, url}
\usepackage[nobreak]{cite}
\usepackage[latin1]{inputenc}
\usepackage[dvips, arrow, matrix, tips, curve]{xy}
\usepackage[british]{babel}
\SelectTips{cm}{10}

\newcommand{\cat}[1]{\mathbf{#1}}

\newcommand{\op}{\mathrm{op}}

\newcommand{\abs}[1]{{\left|{#1}\right|}}

\newcommand{\defeq}{\mathrel{\mathop:}=}

\newcommand{\cd}[2][]{\vcenter{\hbox{\xymatrix#1{#2}}}}

\newcommand{\A}{{\mathcal A}}

\newcommand{\C}{{\mathcal C}}
\newcommand{\D}{{\mathcal D}}
\newcommand{\E}{{\mathcal E}}

\newcommand{\I}{{\mathcal I}}

\newcommand{\K}{{\mathcal K}}

\renewcommand{\P}{{\mathcal P}}

\newcommand{\T}{{\mathcal T}}

\newcommand{\xtor}[1]{\cdl[@1]{{} \ar[r]|-{\object@{|}}^{#1} & {}}}

\makeatletter

\def\hookleftarrowfill@{\arrowfill@\leftarrow\relbar{\relbar\joinrel\rhook}}
\def\twoheadleftarrowfill@{\arrowfill@\twoheadleftarrow\relbar\relbar}
\def\leftbararrowfill@{\arrowdoublefill@{\leftarrow\mkern-5mu}\relbar\mapstochar\relbar\relbar}
\def\Leftbararrowfill@{\arrowdoublefill@{\Leftarrow\mkern-2mu}\Relbar\Mapstochar\Relbar\Relbar}
\def\leftringarrowfill@{\arrowdoublefill@{\leftarrow\mkern-3mu}\relbar{\mkern-3mu\circ\mkern-2mu}\relbar\relbar}
\def\lefttriarrowfill@{\arrowfill@{\mathrel\triangleleft\mkern0.5mu\joinrel\relbar}\relbar\relbar}
\def\Lefttriarrowfill@{\arrowfill@{\mathrel\triangleleft\mkern1mu\joinrel\Relbar}\Relbar\Relbar}

\def\hookrightarrowfill@{\arrowfill@{\lhook\joinrel\relbar}\relbar\rightarrow}
\def\twoheadrightarrowfill@{\arrowfill@\relbar\relbar\twoheadrightarrow}
\def\rightbararrowfill@{\arrowdoublefill@{\relbar\mkern-0.5mu}\relbar\mapstochar\relbar\rightarrow}
\def\Rightbararrowfill@{\arrowdoublefill@{\Relbar\mkern-2mu}\Relbar\Mapstochar\Relbar\Rightarrow}
\def\rightringarrowfill@{\arrowdoublefill@\relbar\relbar{\mkern-2mu\circ\mkern-3mu}\relbar{\mkern-3mu\rightarrow}}
\def\righttriarrowfill@{\arrowfill@\relbar\relbar{\relbar\joinrel\mkern0.5mu\mathrel\triangleright}}
\def\Righttriarrowfill@{\arrowfill@\Relbar\Relbar{\Relbar\joinrel\mkern1mu\mathrel\triangleright}}

\def\leftrightarrowfill@{\arrowfill@\leftarrow\relbar\rightarrow}
\def\mapstofill@{\arrowfill@{\mapstochar\relbar}\relbar\rightarrow}

\renewcommand*\xleftarrow[2][]{\ext@arrow 20{20}0\leftarrowfill@{#1}{#2}}
\providecommand*\xLeftarrow[2][]{\ext@arrow 60{22}0{\Leftarrowfill@}{#1}{#2}}
\providecommand*\xhookleftarrow[2][]{\ext@arrow 10{20}0\hookleftarrowfill@{#1}{#2}}
\providecommand*\xtwoheadleftarrow[2][]{\ext@arrow 60{20}0\twoheadleftarrowfill@{#1}{#2}}
\providecommand*\xleftbararrow[2][]{\ext@arrow 10{22}0\leftbararrowfill@{#1}{#2}}
\providecommand*\xLeftbararrow[2][]{\ext@arrow 50{24}0\Leftbararrowfill@{#1}{#2}}
\providecommand*\xleftringarrow[2][]{\ext@arrow 10{26}0\leftringarrowfill@{#1}{#2}}
\providecommand*\xlefttriarrow[2][]{\ext@arrow 80{24}0\lefttriarrowfill@{#1}{#2}}
\providecommand*\xLefttriarrow[2][]{\ext@arrow 80{24}0\Lefttriarrowfill@{#1}{#2}}

\renewcommand*\xrightarrow[2][]{\ext@arrow 01{20}0\rightarrowfill@{#1}{#2}}
\providecommand*\xRightarrow[2][]{\ext@arrow 04{22}0{\Rightarrowfill@}{#1}{#2}}
\providecommand*\xhookrightarrow[2][]{\ext@arrow 00{20}0\hookrightarrowfill@{#1}{#2}}
\providecommand*\xtwoheadrightarrow[2][]{\ext@arrow 03{20}0\twoheadrightarrowfill@{#1}{#2}}
\providecommand*\xrightbararrow[2][]{\ext@arrow 01{22}0\rightbararrowfill@{#1}{#2}}
\providecommand*\xRightbararrow[2][]{\ext@arrow 04{24}0\Rightbararrowfill@{#1}{#2}}
\providecommand*\xrightringarrow[2][]{\ext@arrow 01{26}0\rightringarrowfill@{#1}{#2}}
\providecommand*\xrighttriarrow[2][]{\ext@arrow 07{24}0\righttriarrowfill@{#1}{#2}}
\providecommand*\xRighttriarrow[2][]{\ext@arrow 07{24}0\Righttriarrowfill@{#1}{#2}}

\providecommand*\xmapsto[2][]{\ext@arrow 01{20}0\mapstofill@{#1}{#2}}
\providecommand*\xleftrightarrow[2][]{\ext@arrow 10{22}0\leftrightarrowfill@{#1}{#2}}
\providecommand*\xLeftrightarrow[2][]{\ext@arrow 10{27}0{\Leftrightarrowfill@}{#1}{#2}}

\makeatother

\newcommand{\twocong}[2][0.5]{\ar@{}[#2] \save ?(#1)*{\cong}\restore}
\newcommand{\twoeq}[2][0.5]{\ar@{}[#2] \save ?(#1)*{=}\restore}
\newcommand{\rtwocell}[3][0.5]{\ar@{}[#2] \ar@{=>}?(#1)+/l 0.2cm/;?(#1)+/r 0.2cm/^{#3}}
\newcommand{\ltwocell}[3][0.5]{\ar@{}[#2] \ar@{=>}?(#1)+/r 0.2cm/;?(#1)+/l 0.2cm/^{#3}}
\newcommand{\ltwocello}[3][0.5]{\ar@{}[#2] \ar@{=>}?(#1)+/r 0.2cm/;?(#1)+/l 0.2cm/_{#3}}
\newcommand{\dtwocell}[3][0.5]{\ar@{}[#2] \ar@{=>}?(#1)+/u  0.2cm/;?(#1)+/d 0.2cm/^{#3}}
\newcommand{\dltwocell}[3][0.5]{\ar@{}[#2] \ar@{=>}?(#1)+/ur  0.2cm/;?(#1)+/dl 0.2cm/^{#3}}
\newcommand{\drtwocell}[3][0.5]{\ar@{}[#2] \ar@{=>}?(#1)+/ul  0.2cm/;?(#1)+/dr 0.2cm/^{#3}}
\newcommand{\dthreecell}[3][0.5]{\ar@{}[#2] \ar@3{->}?(#1)+/u  0.2cm/;?(#1)+/d 0.2cm/^{#3}}
\newcommand{\utwocell}[3][0.5]{\ar@{}[#2] \ar@{=>}?(#1)+/d 0.2cm/;?(#1)+/u 0.2cm/_{#3}}
\newcommand{\dtwocelltarg}[3][0.5]{\ar@{}#2 \ar@{=>}?(#1)+/u  0.2cm/;?(#1)+/d 0.2cm/^{#3}}
\newcommand{\utwocelltarg}[3][0.5]{\ar@{}#2 \ar@{=>}?(#1)+/d  0.2cm/;?(#1)+/u 0.2cm/_{#3}}

\newdir{(}{{}*!<0em,-.14em>-\cir<.14em>{l^r}}
\newdir{ (}{{}*!/-5pt/\dir{(}}
\newdir{ >}{{}*!/-5pt/\dir{>}}

\swapnumbers
\theoremstyle{plain}
\newtheorem{Thm}{Theorem}[section]
\newtheorem{Prop}[Thm]{Proposition}

\theoremstyle{definition}

\theoremstyle{remark}

\DeclareMathOperator{\Lan}{Lan}

\begin{document}
 \leftmargini=2em
\title{A characterisation of algebraic exactness}
\author{Richard Garner}
\address{Department of Computing, Macquarie University, NSW 2109, Australia}
\email{richard.garner@mq.edu.au} \subjclass[2010]{}
\date{\today}
\begin{abstract}
An \emph{algebraically exact} category is one that admits all of the limits and
colimits which every variety of algebras possesses and every forgetful functor between varieties preserves, and which verifies the same interactions between these limits and colimits as hold in any variety. Such categories were studied by Ad{\'a}mek, Lawvere and Rosick{\'y}: they characterised them as
the categories with small limits and sifted colimits for which the functor
taking sifted colimits is continuous. They conjectured that a complete and
sifted-cocomplete category should be algebraically exact just when it is
Barr-exact, finite limits commute with filtered colimits, regular epimorphisms
are stable by small products, and filtered colimits distribute over small
products. We prove this conjecture.
\end{abstract}
\thanks{The support of the Australian Research Council and DETYA is gratefully acknowledged}
 \maketitle

\newcommand{\sind}{\mathcal S \mathrm{ind}}
\newcommand{\sindk}{\mathcal S_{\kappa'}}

\section{Introduction}

In the series of papers~\cite{Adamek2001How-algebraic,Adamek2001On-algebraically,Adamek2004Toward} was introduced and studied the notion of an
\emph{algebraically exact} category. A category $\C$ is said to be
algebraically exact if, firstly, it admits all of the operations $\C^\A \to \C$ of small arity which every variety of (finitary, many-sorted) algebras supports and every forgetful functor between varieties preserves, and secondly, it obeys all of the equations
between such operations as are satisfied in every variety. Any variety admits
small limits and sifted colimits, and every forgetful functor between varieties preserves them; recall from~\cite{Adamek2001On-sifted} that \emph{sifted colimits} are those
which commute with finite products in $\cat{Set}$, most important amongst these
being the filtered colimits, and the coequalisers of reflexive pairs. It
follows that any algebraically exact category also admits small limits and
sifted colimits; and it turns out that these two kinds of operations in fact
generate all of those required of an algebraically exact category. As regarding
the equations that hold between these operations, we observe that in any
variety, the following four exactness properties are verified:
\begin{enumerate}[(E1)]
\item Regular epimorphisms are stable under pullback, and equivalence relations are effective (i.e., the category is Barr-exact);
\item Finite limits commute with filtered colimits;
\item Regular epimorphisms are stable by small products;
\item Filtered colimits distribute over small products.
\end{enumerate}

It follows that these same conditions are verified in any algebraically exact
category, and it was conjectured in~\cite{Adamek2001How-algebraic} that, in
fact, these four conditions completely characterise the algebraically exact
categories amongst those categories with small limits and sifted colimits. The
conjecture was proved in~\cite{Adamek2001On-algebraically} for the case of
cocomplete categories with a regular generator, and in~\cite{Adamek2004Toward}
for the case of categories with finite coproducts; the purpose of this article
is to prove it in its full generality. We shall do so using techniques
developed in~\cite{Garner2011Lex-colimits}, though the arguments are
straightforward enough that we can reproduce them in full here, so making this
article entirely self-contained.

In order to state the conjecture more precisely, we will make use of a different description of the algebraically exact categories. We recall from~\cite{Adamek2001On-sifted} the construction which to every locally small category $\C$ assigns its free completion $\sind(\C)$ under sifted colimits. As in~\cite[Theorem 5.35]{Kelly1982Basic}, we may obtain $\sind(\C)$ as the closure of the representables in $[\C^\op, \cat{Set}]$ under sifted colimits, and now the restricted Yoneda embedding $W \colon \C \to \sind(\C)$ provides the unit at $\C$ of a Kock-Z\"oberlein pseudomonad~\cite{Kock1995Monads} on $\cat{CAT}$, whose pseudoalgebras are the sifted-cocomplete categories. Thus a category $\C$ admits sifted colimits just when $W \colon \C \to \sind(\C)$ admits a left adjoint.

It was shown in~\cite[Theorem 3.11]{Adamek2001How-algebraic} that if $\C$ is complete, then so too is $\sind(\C)$; that if $F \colon \C \to \D$ is a continuous functor between complete categories, then so too is $\sind(F)$; and that the unit $\C \to \sind(\C)$ and multiplication $\sind(\sind(\C)) \to \sind(\C)$ are always continuous functors. It follows that the  pseudomonad $\sind$ restricts and corestricts to one on $\cat{CONTS}$, the $2$-category of complete categories and continuous functors; and it was shown in~\cite[Corollary 4.4]{Adamek2001How-algebraic} that the pseudoalgebras for this restricted pseudomonad are precisely the algebraically exact categories described above. Thus a complete and sifted-cocomplete category $\C$ is algebraically exact just when $W \colon \C \to \sind(\C)$ admits a left adjoint which is \emph{continuous}.
For the purposes of this paper, we will take this last as our definition of an algebraically exact category; and our goal, then, is to prove:
\begin{Thm}\label{totalthm}
A complete and sifted-cocomplete category $\C$ is algebraically exact just when it satisfies conditions (E1)--(E4).
\end{Thm}
In fact, as remarked above, any algebraically exact category does indeed
satisfy (E1)--(E4); and so our task is to show that these conditions in turn
imply algebraic exactness.


\section{The result}
The basic idea behind the proof of Theorem~\ref{totalthm} is to show that any
category $\C$ satisfying (E1)--(E4) admits a full structure-preserving
embedding into some $\E$ which is an essential localisation of a presheaf
topos. Any such $\E$ will be algebraically exact; and now we may reflect this
property  along the full embedding, so concluding that $\C$ itself is
algebraically exact. This argument does not quite work as it stands, for
reasons of size. The $\E$ into which we would like to embed is a topos of
sheaves on $\C$, but only when $\C$ is small may such a topos be constructed;
in which situation, with $\C$ being small, and also small-complete, it is
necessarily a preorder, which is far too restrictive. To overcome this problem,
we will first prove a variant of Theorem~\ref{totalthm}, in which suitable
bounds have been introduced on the size of the limits and colimits required,
and then deduce the general result from this.

Our cardinality bounds will be governed by an infinite regular cardinal $\kappa$. Given
any such $\kappa$, we define $\kappa'$ to be the cardinal $(\Sigma_{\gamma < \kappa} 2^{\gamma})^+$,
and the pair $(\kappa, \kappa')$ now has the property that whenever $\mu <
\kappa$ and $\lambda < \kappa'$, we have $\lambda^\mu < \kappa'$:
see~\cite[Proposition 2.3.5]{Makkai1989Accessible}. By a $\kappa$-limit we
shall mean one indexed by a diagram of cardinality $< \kappa$, and we attach a
corresponding meaning to the term $\kappa'$-colimit.  We shall now describe a
variant of the notion of algebraic exactness, which we term
\emph{$\kappa$-algebraic exactness}, that deals only with $\kappa$-limits and
$\kappa'$-colimits.

\looseness=-1 There is a slight delicacy here as to the kinds of
$\kappa'$-colimit we will consider. The obvious choice would be the sifted
$\kappa'$-colimits---which we emphasise means the $\kappa'$-small sifted
colimits, and \emph{not} the colimits which commute in $\cat{Set}$ with $\kappa'$-small
products---but this choice is in fact inappropriate. It follows from~\cite[Proposition 5.1]{Adamek2001On-algebraically} that if $\C$ is complete
then $\sind(\C)$ is the closure of the representables in $[\C^\op, \cat{Set}]$ under
reflexive coequalisers and filtered colimits, so that a complete $\C$ admits
sifted colimits just when it admits reflexive coequalisers and filtered
colimits. When we bound the cardinality of our colimits, it turns out to be the
reflexive coequalisers together with the filtered $\kappa'$-colimits which are relevant,
and not the sifted $\kappa'$-colimits; recall from~\cite{Adamek2010What} that
the latter class of colimits is in general \emph{strictly} larger.

We consider the $2$-category $\kappa\text-\cat{CONTS}$ of $\kappa$-complete
categories and $\kappa$-continuous functors between them; on this, we will
describe a pseudomonad whose pseudoalgebras will be the $\kappa$-algebraically
exact categories we seek to define. Observe first that as well as the
pseudomonad $\sind$ on $\cat{CAT}$ we also have the pseudomonad $\P$ which
freely adds small colimits. Proposition 4.3 and Remark 6.6
of~\cite{Day2007Limits} prove that if $\C$ is $\kappa$-complete, then so is $\P
\C$; that if $F \colon \C \to \D$ is a $\kappa$-continuous functor between such
categories, then so is $\P F$; and that $\P$'s unit and multiplication are
always $\kappa$-continuous. Thus we may restrict and corestrict $\P$ to a
pseudomonad on $\kappa\text-\cat{CONTS}$; and the pseudomonad of interest to us
will be a submonad of this, defined as follows. For each $\C$ in
$\kappa\text-\cat{CONTS}$, we let $\sindk(\C)$ denote the closure of $\C$ in
$\P\C$ under $\kappa$-limits, reflexive coequalisers, and filtered
$\kappa'$-colimits, and let $V \colon \C \to \sindk(\C)$ denote the restricted
Yoneda embedding. Now~\cite[Proposition 3.1]{Garner2011Lex-colimits} ensures
that this $V$ provides the unit at $\C$ of a Kock-Z\"oberlein pseudomonad on
$\kappa\text-\cat{CONTS}$; and a $\kappa$-algebraically exact category will be,
by definition, a pseudoalgebra for this pseudomonad. In other words, a
$\kappa$-complete category $\C$ is \emph{$\kappa$-algebraically exact} just
when the embedding $V \colon \C \to \sindk(\C)$ admits a $\kappa$-continuous
left adjoint. Observe that this implies that $\C$ has reflexive coequalisers
and filtered $\kappa'$-colimits, but may not imply that it has all sifted
$\kappa'$-colimits; this is in accordance with the remarks of the preceding
paragraph.

We shall now prove the following refinement of Theorem~\ref{totalthm}.

\begin{Thm}\label{mainthm}
A category $\C$ with $\kappa$-limits, reflexive coequalisers and filtered
$\kappa'$-colimits is $\kappa$-algebraically exact just when:
\begin{enumerate}[(E1')]
\item It is Barr-exact;
\item Finite limits commute with filtered $\kappa'$-colimits;
\item Regular epimorphisms are stable by $\kappa$-small products;
\item Filtered $\kappa'$-colimits distribute over $\kappa$-small products.
\end{enumerate}
\end{Thm}

Clearly, a complete and sifted-cocomplete $\C$ satisfies (E1')--(E4') for each
regular $\kappa$ if and only if it satisfies (E1)--(E4). On the other hand, we
have:
\begin{Prop}\label{propreduce}
A complete and sifted-cocomplete category $\C$ is algebraically exact if and only
if it is $\kappa$-algebraically exact for each regular $\kappa$.
\end{Prop}
By virtue of this Proposition and the comment preceding it, we may prove Theorem~\ref{totalthm} by proving Theorem~\ref{mainthm}, and then taking the conjunction of all its instances as $\kappa$ ranges across the small regular cardinals.
\begin{proof}[Proof of Proposition~\ref{propreduce}]
For every $\kappa$, we observe that $\sind(\C)$ is closed under $\kappa$-limits, reflexive coequalisers and filtered $\kappa'$-colimits in $[\C^\op,
\cat{Set}]$; whence $\sindk(\C) \subset \sind(\C)$ with the inclusion
preserving all $\kappa$-limits. Hence if $W \colon \C \to \sind(\C)$ admits a
continuous left adjoint, then by restriction each $V \colon \C \to \sindk(\C)$
will admit a $\kappa$-continuous left adjoint.

Conversely, suppose that each $V \colon \C \to \sindk(\C)$ admits a
$\kappa$-continuous left adjoint. As observed above, since $\C$ is complete, it
follows by~\cite[Proposition 5.1]{Adamek2001On-algebraically} that $\sind(\C)$
is the closure of the representables in $[\C^\op, \cat{Set}]$ under reflexive
coequalisers and filtered colimits. But it is easy to see that the collection
of $\varphi \in \sind(\C)$ which lie in some $\sindk(\C)$ contains the
representables and is closed under reflexive coequalisers and filtered
colimits, and so must be all of $\sind(\C)$; which is to say that $\sind(\C) =
\bigcup_\kappa \sindk(\C)$. Thus, since each $V \colon \C \to \sindk(\C)$
admits a left adjoint, so too does $W \colon \C \to \sind(\C)$, and it remains
to show that this left adjoint is continuous. Given a small diagram $D \colon
\I \to \sind(\C)$, we may choose a regular cardinal $\kappa$ such that
 $DI \in \sindk(\C)$ for each $I \in \I$ and also $\abs{\I} < \kappa$; now
the diagram $D$ factors as $D' \colon \I \to \sindk(\C)$, and the left adjoint
of $\C \to \sindk(\C)$ preserves the limit of $D'$: from which it follows that
the left adjoint of $W$ preserves that of $D$, as required.
\end{proof}

We now prove Theorem~\ref{mainthm} for the case of a small $\C$. Given such a
$\C$ satisfying the conditions of the theorem, we shall embed it into a $\kappa$-algebraically exact category via a functor preserving $\kappa$-limits, reflexive
coequalisers and filtered $\kappa'$-colimits.
It will then follow that $\C$ is $\kappa$-algebraically exact by virtue of the following result. 
\begin{Prop}\label{embedding}
Let $J \colon \C \to \E$ be fully faithful; suppose moreover that $\C$ has, and that $J$ preserves, $\kappa$-limits, reflexive coequalisers and filtered $\kappa'$-colimits,
and that $\E$ is $\kappa$-algebraically exact. Then $\C$ is also $\kappa$-algebraically exact.
\end{Prop}
\begin{proof}
Because $\E$ is $\kappa$-algebraically exact, the functor $J$ admits a left Kan extension
\begin{equation*}\cd[@-0.5em]{
 \C \ar[d]_V \ar[r]^J \twocong[0.3]{dr}{} & \E\\
 \sindk(\C) \ar[ur]_{\Lan_V J} & {}
 }
\end{equation*}
 along $V$, which may be calculated as the composite
\begin{equation*}
\sindk(\C) \xrightarrow{\sindk(J)} \sindk(\E) \xrightarrow{\quad L \quad} \E
\end{equation*}
with $L$ the $\kappa$-continuous left adjoint of $V \colon \E \to \sindk(\E)$.
Now $\sindk(J)$ is an algebra morphism between free $\sindk$-algebras, and as
such, preserves $\kappa$-limits, reflexive coequalisers and filtered $\kappa'$-colimits; whilst $L$ preserves all colimits,
being a left adjoint. It follows that $\Lan_V J$, like $J$, preserves $\kappa$-limits, reflexive coequalisers and filtered $\kappa'$-colimits; whence the collection of
$\varphi \in \sindk(\C)$ for which $\Lan_V J$ lands in the essential image of
$J$ contains the representables and is closed under $\kappa$-limits, reflexive coequalisers and filtered $\kappa'$-colimits, and so must be all of $\sindk(\C)$. Hence $\Lan_V J$
factors through $J$, up-to-isomorphism; and the factorisation $\sindk(\C) \to
\C$ so induced, which is clearly $\kappa$-continuous, may also be shown to be
left adjoint to $V \colon \C \to \sindk(\C)$, so that $\C$ is indeed
$\kappa$-algebraically exact.
\end{proof}
Given a small, $\kappa$-complete $\C$, admitting reflexive coequalisers and
filtered $\kappa'$-colimits, and satisfying (E1')--(E4'), we now exhibit an
embedding of the above form; as anticipated at the start of this section, it will in fact be an embedding into a topos. We consider the smallest topology on $\C$ for which all regular epimorphisms are covering, and for which the
colimit injections into each filtered $\kappa'$-colimit are covering.
(E1') and (E2') ensure that this topology is subcanonical and
so we have a full embedding $J \colon \C \to \cat{Sh}(\C)$.
\begin{Prop}
The full embedding $J \colon \C \to \cat{Sh}(\C)$ preserves $\kappa$-limits, reflexive coequalisers and filtered $\kappa'$-colimits.
\end{Prop}
\begin{proof}
Clearly $J$ preserves all limits that exist, so in particular $\kappa$-limits. It also preserves regular epimorphisms, since the given topology contains the regular one, and we will show below that it preserves filtered $\kappa'$-colimits. It will then follow that it preserves reflexive coequalisers too, since in $\C$ and in $\cat{Sh}(\C)$, we may exploit (E1') and (E2') to construct such coequalisers from finite limits, countable filtered colimits and coequalisers of equivalence relations, all of which are preserved by $J$; the argument is standard and given in precisely the form we need in~\cite[Theorem
2.6]{Adamek2004Toward}.

It remains to show that $J$ preserves filtered $\kappa'$-colimits. Observe that if $(p_k \colon Dk \to X \mid k \in \K)$ is such a colimit in $\C$, then $J$ will preserve it just when every sheaf $\C^\op \to \cat{Set}$  sends it to a limit in $\cat{Set}$. Let $F$ be such a sheaf. Since the family $(p_k \mid k \in \K)$ is covering, we may identify $FX$ with the set of matching families for this covering. In other words, if
\begin{equation*}
\cd[@-0.5em]{
  D_{jk} \ar[r]^-{d_{jk}} \ar[d]_{c_{jk}} & Dj \ar[d]^{p_j} \\
  Dk \ar[r]_{p_k} & X
}  
\end{equation*}
is a pullback for each $j, k \in \K$, then we may identify $FX$ with the set
\begin{equation}\label{eq:theset}\tag{$\ast$}
 \{ \vec x \in \Pi_k FDk \,\mid\,  Fd_{jk}(x_j) = Fc_{jk}(x_k) \text{ for all $j, k \in \K$}\}\rlap{ .}
\end{equation}
Under this identification, the canonical comparison map $FX \to \lim FD$ is just the inclusion between these sets, seen as subobjects of $\Pi_k FDk$, and so injective; it remains to show that it is also surjective. Thus we must show that each $\vec x \in \lim FD$ lies in~\eqref{eq:theset}, or in other words, that  $Fd_{jk}(x_j) = Fc_{jk}(x_k)$ for each such $\vec x$ and each $j, k \in J$. 
To this end, we consider the category $\K'$ of cospans from $j$ to $k$ in $\K$; since $\K$ is filtered and $\kappa'$-small, it follows easily that $\K'$ is too. We define a functor $E \colon \K' \to \C$ by sending each cospan $f \colon j \to \ell \leftarrow k \colon g$ in $\K'$ to the 
%
%
apex of the pullback square\[
\cd[@-0.5em]{
  E(f,g) \ar[r]^-{u_{f,g}} \ar[d]_{v_{f,g}} & Dj \ar[d]^{Df} \\
  Dk \ar[r]_{Dg} & D\ell
}
\]
in $\C$. A simple calculation shows that
$p_k.v_{f,g} = p_j.u_{f,g}$, so that we have induced maps $q_{f,g} \defeq (u_{f,g}, v_{f,g}) \colon E(f,g) \to D_{jk}$, constituting a cocone $q$ under $E$ with vertex $D_{jk}$. We claim that this cocone is colimiting; whereupon, 
by the preceding part of the argument, the comparison $FD_{jk} \to \lim FE$ induced by $q$ will be  monic, and so the family $(Fq_{f,g} \mid (f,g) \in \K')$ jointly monic. Thus in order to verify that $Fd_{jk}(x_j) = Fc_{jk}(x_k)$, and so complete the proof, it will be enough to observe that for each $f \colon j \to \ell \leftarrow k \colon g$ in $\K'$, we have:
\begin{align*}Fq_{f,g}(Fd_{jk}(x_j)) &= Fu_{f,g}(x_j) = Fu_{f,g}(FDf(x_\ell))\\ & = Fv_{f,g}(FDg(x_\ell)) = Fv_{f,g}(x_k) \\ &= Fq_{f,g}(Fc_{jk}(x_k))\rlap{ .}
\end{align*}

It remains to verify that $q$ is colimiting. For this, let $V \colon \K' \to \K$ denote the functor sending a $j,k$-cospan to its central object, and $\iota_1 \colon \Delta j \to V \leftarrow \Delta k \colon \iota_2$ the evident natural transformations. Now we have a commutative cube
\begin{equation*}
\cd[@!@-2.5em@C+0.5em]{
  E \ar[rr]^u \ar[dd]_v \ar[dr]^{q} & & \Delta Dj \ar[dd]_(0.25){D\iota_1} \ar@{=}[dr] \\ &
  \Delta(D_{jk}) \ar[rr]_(0.7){\Delta d_{jk}} \ar[dd]^(0.75){\Delta c_{jk}} & &
  \Delta(Dj) \ar[dd]^{\Delta p_j} \\
  \Delta(Dk) \ar[rr]^(0.33){D \iota_2} \ar@{=}[dr] & & DV \ar[dr]^{pV} \\ &
  \Delta(Dk) \ar[rr]_{\Delta p_k} & & \Delta X
}
\end{equation*}
in $[\K', \C]$; its front and rear faces are pullbacks, and by (E2') will remain so on applying the functor $\mathrm{colim} \colon [\K', \C] \to \C$. To show that $q$ is colimiting is equally to show that it is inverted by $\mathrm{colim}$; for which, by the previous sentence, it is enough to show that $pV$ is likewise inverted. But $\K$'s filteredness implies easily that $V \colon \K' \to \K$ is a final functor, so that $pV$, like $p$, is a colimiting cocone, and so inverted by $\mathrm{colim}$ as required.
\end{proof}
We thus have a full structure-preserving embedding $\C \to \cat{Sh}(\C)$ and
the only thing left to verify is that $\cat{Sh}(\C)$ is in fact $\kappa$-algebraically exact.
The key to doing so is the following proposition.
\begin{Prop}
If $\E$ is reflective in a presheaf category via a $\kappa$-continuous reflector, then $\E$ is $\kappa$-algebraically exact.
\end{Prop}
\begin{proof}
If $\C$ is small, then $\P \C = [\C^\op, \cat{Set}]$, and now the restricted
Yoneda embedding $\P \C \to \P \P \C$ admits a continuous left adjoint $\P \P
\C \to \P \C$, this being the multiplication at $\C$ of the pseudomonad $\P$.
Since $\sindk(\P \C)$ is closed in $\P \P \C$ under $\kappa$-limits, it follows
by restriction that $\P \C \to \sindk(\P \C)$ admits a $\kappa$-continuous left
adjoint; and so every presheaf category is $\kappa$-algebraically exact. Now if
$\E$ is reflective in the $\kappa$-algebraically exact $[\C^\op, \cat{Set}]$
via a $\kappa$-continuous reflector, then it is an adjoint retract of $[\C^\op,
\cat{Set}]$ in $\kappa\text-\cat{CONTS}$, and so by a standard property of
Kock-Z\"oberlein pseudomonads, must itself be $\kappa$-algebraically exact.
\end{proof}
Thus it is enough to show that $\cat{Sh}(\C)$ is reflective in
$[\C^\op, \cat{Set}]$ via a $\kappa$-continuous reflector. This will be a
consequence of the following result, which may be found proven---though with ``small'' harmlessly replacing our ``$\kappa$-small''---in~\cite[Theorem 4.2]{Kelly1989On-the-complete}; we shall not recall the details, since we shall not need them in what follows.
\begin{Prop}
A left exact reflector $L \colon [\C^\op, \cat{Set}] \to \E$ preserves all $\kappa$-small limits if and only if the covering sieves for the corresponding topology are closed under $\kappa$-small intersections in $[\C^\op, \cat{Set}]$.
\end{Prop}
We are therefore required to show that any $\kappa$-small intersection of covering sieves for the above-defined topology on $\C$ is again covering. Clearly it is sufficient to consider the case where the sieves participating in the intersection are generating ones for the topology. We can decompose any such intersection of sieves as an intersection
\begin{equation*}
\bigcap_{i \in I} \mathcal S_i \cap \bigcap_{j \in J} \T_j
\end{equation*}
where each indexing set $I$ and $J$ is $\kappa$-small, each sieve $\mathcal S_i$ is generated by a regular epimorphism $e_i \colon A_i \twoheadrightarrow X$ and each sieve $\T_j$ is generated by a $\kappa'$-small filtered colimit cocone $((q_j)_{x} \colon D_j(x) \to X \mid x \in \A_j)$.

Now we can form the $\kappa$-small product $\Pi_i e_i \colon \Pi_i A_i \to
\Pi_i X$; by condition (E3') this is a regular epimorphism in $\C$, and  by
regularity, so also is its pullback $e \colon A \to X$ along the diagonal $X
\to \Pi_i X$. Clearly a map $Z \to X$ factors through $e$ just when it factors
through each $e_i$, and so the covering sieve $\mathcal S$ generated by $e$ is
the intersection $\bigcap_i \mathcal S_i$.

In a similar manner, we can form the filtered category $\Pi_j \A_j$; since
$\abs{J} < \kappa$, and each $\abs{\A_j} < \kappa'$, we have also that
$\abs{\Pi_j \A_j} < \kappa'$. Now on considering the diagram $D \colon \Pi_j
\A_j \to \C$ defined by $D(x_j \mid j \in J) = \Pi_j D_j(x_j)$, condition (E4')
asserts that $\Pi_j X$ is a colimit for it; so that on pulling back along the
diagonal $X \to \Pi_j X_j$, we conclude that $X$ is a colimit for the diagram
$D' \colon \Pi_j \A_j \to \C$ which sends $(x_j \mid j \in J)$ to the fibre
product of the maps $(q_j)_{x_j} \colon D_j(x_j) \to X$. Now we see as before
that the covering sieve $\T$ generated by this filtered $\kappa'$-colimit
cocone is precisely $\bigcap_j \T_j$.

It follows that $\bigcap_{i} \mathcal S_i \cap \bigcap_{j} \T_j = \mathcal S
\cap \T$ is a covering sieve, since covering sieves are always closed under
finite intersections, and this  completes the proof of:
\begin{Prop}
If the small, $\kappa$-complete $\C$ with reflexive coequalisers and filtered
$\kappa'$-colimits satisfies (E1')--(E4'), then it admits a full
structure-preserving embedding into a $\kappa$-algebraically exact category,
and so is itself $\kappa$-algebraically exact.
\end{Prop}

It remains to prove Theorem~\ref{mainthm} for categories of no matter what
size. So let $\C$ be a category with
$\kappa$-limits, reflexive coequalisers and filtered $\kappa'$-colimits, satisfying (E1')--(E4').
We call a full, replete subcategory \emph{$\kappa$-closed} if it is closed in $\C$ under
the limits and colimits just mentioned. Clearly, each small, $\kappa$-closed subcategory of $\C$
satisfies (E1')--(E4'), and so by the preceding proposition is $\kappa$-algebraically exact.
We may now conclude that the same is true of $\C$ by way of the following result.

\begin{Prop}
A $\kappa$-complete  $\C$ admitting reflexive coequalisers and filtered
$\kappa'$-colimits is $\kappa$-algebraically exact so long as all of its small
$\kappa$-closed subcategories are.
\end{Prop}
\begin{proof}
Suppose that each
$\kappa$-closed subcategory of $\C$ is $\kappa$-algebraically exact; we must show
that $\C$ is too, or in other words, that $V \colon \C \to \sindk(\C)$ admits a
$\kappa$-continuous left adjoint.
To this end, consider the collection of $\varphi \in \sindk(\C)$ for which
there exists a small $\kappa$-closed $J \colon \D \hookrightarrow \C$ with
$\varphi$ lying in the essential image of the fully faithful $\sindk(J) \colon
\sindk(\D) \to \sindk(\C)$. It is easy to show that this collection
contains the representables and is closed under $\kappa$-limits, reflexive coequalisers and filtered
$\kappa'$-colimits, and so is all of $\sindk(\C)$. It follows that $\C \to
\sindk(\C)$ admits a left adjoint, since each $\D \to \sindk(\D)$ does by
assumption.

To show that this left adjoint is moreover $\kappa$-continuous, consider a
$\kappa$-small diagram $X \colon \I \to \sindk(\C)$. For each $I \in \I$ we can
find a small $\kappa$-closed $\D_I \subset \C$ with $XI$ in the essential image
of $\sindk(\D_I) \to \sindk(\C)$; now taking $\D$ to be the closure
of $\bigcup_I \D_I$ in $\C$ under $\kappa$-limits, reflexive coequalisers and
filtered $\kappa'$-colimits, we obtain another small $\kappa$-closed
subcategory. The diagram $X$ factors up-to-isomorphism through the fully faithful
$\sindk(\D) \rightarrow \sindk(\C)$ as $X' \colon \I \to \sindk(\D)$, say; and
now by assumption, the left adjoint of $\D \to \sindk(\D)$ preserves the limit
of $X'$, whence the left adjoint of $\C \to \sindk(\C)$ preserves that of $X$,
as required.
\end{proof}
	
This completes the proof of Theorem~\ref{mainthm} for categories of any size; and now,
as discussed previously, taking the conjunction of all instances of this theorem as $\kappa$ ranges over
the small regular cardinals completes the proof of Theorem~\ref{totalthm}.




\bibliographystyle{acm}

\bibliography{bibdata}

\end{document}